\documentclass[12pt,reqno]{amsart}

\usepackage{amssymb}

\usepackage{eucal}

\input{diagrams}

\def\cC{{\mathcal C}}

\def\fm{{\mathfrak m}}

\def\sD{\mbox{\sf D}}

\def\CIdim{\operatorname{CI-dim}}

\def\Coker{\operatorname{Coker}}
\def\colim{\operatorname{colim}}
\def\D{\sD}

\def\dim{\operatorname{dim}}

\def\Ext{\operatorname{Ext}}

\def\Gammam{\Gamma_{\fm}}
\def\Gammamopp{\Gamma_{{\fm}^{\opp}}}

\def\H{\operatorname{H}}

\def\Hom{\operatorname{Hom}}
\def\id{\operatorname{id}}

\def\lim{\operatorname{lim}}

\def\LTensor{\stackrel{\operatorname{L}}{\otimes}}

\def\opp{\operatorname{op}}

\def\R{\operatorname{R}}

\def\RGammam{\operatorname{R}\!\Gamma_{\fm}}
\def\RGammamopp{\operatorname{R}\!\Gamma_{{\fm}^{\opp}}}

\def\RHom{\operatorname{RHom}}

\def\Tor{\operatorname{Tor}}

\numberwithin{equation}{part}

\newtheorem{Lemma}{Lemma}[section]
\newtheorem{Theorem}[Lemma]{Theorem}

\newtheorem{Proposition}[Lemma]{Proposition}

\theoremstyle{definition}
\newtheorem{Definition}[Lemma]{Definition}
\newtheorem{Setup}[Lemma]{Setup}

\newtheorem{Remark}[Lemma]{Remark}

\begin{document}

\title[Symmetry for Ext vanishing]
{Symmetry theorems for Ext vanishing}

\author{Peter J\o rgensen}
\address{Department of Pure Mathematics\\
         University of Leeds\\
         Leeds LS2 9JT\\
         United Kingdom}
\email{popjoerg@maths.leeds.ac.uk}
\urladdr{http://www.maths.leeds.ac.uk/\~{ }popjoerg}




\begin{abstract}

It was proved by Avramov and Buchweitz that if $A$ is a commutative
local complete intersection ring with finitely ge\-ne\-ra\-ted modules
$M$ and $N$, then the $\Ext$ groups between $M$ and $N$ vanish from
some step if and only if the $\Ext$ groups between $N$ and $M$ vanish
from some step.

This paper shows that the same is true under the weaker conditions
that $A$ is Gorenstein and that $M$ and $N$ have finite complete
intersection dimension.  The result is also proved if $A$ is
Gorenstein and has finite Cohen-Macaulay type.

Similar results are given for two types of non-commutative rings:
Frobenius algebras and complete semi-local algebras.

\end{abstract}

\maketitle

\setcounter{section}{-1}
\section{Introduction}
\label{sec:introduction}

Let $A$ be a commutative local complete intersection ring with
finitely generated modules $M$ and $N$.  It is a surprising result of
\cite{AB} that symmetry of Ext vanishing holds in the sense that
\begin{equation}
\label{equ:symmetry}
  \Ext_A^i(M,N) = 0 \; \mbox{for} \; i \gg 0
  \; \Leftrightarrow \;
  \Ext_A^i(N,M) = 0 \; \mbox{for} \; i \gg 0.
\end{equation}
This paper proves that \eqref{equ:symmetry} remains true if $A$ is a
commutative local Gorenstein ring and $M$ and $N$ have finite complete
intersection dimension in the sense of \cite{AGP}.  

These conditions are weaker than the ones in \cite{AB} because a
complete intersection ring is a Gorenstein ring for which each
finitely generated module has finite complete intersection dimension.

It is also proved that \eqref{equ:symmetry} is true if $A$ is a
commutative local Gorenstein ring of finite Cohen-Macaulay type.

The method of the paper is to isolate two simple homological
properties of (bi)\-mo\-du\-les, property (R) (``Rigidity'') and
property (S) (``Symmetry''), which make possible the abstract $\Ext$
vanishing result theorem \ref{thm:main}.  This in turn implies the
results already stated, see theorem \ref{thm:commutative}, and also
two results dealing with non-commutative rings: Theorem
\ref{thm:Frobenius} on certain Frobenius algebras and theorem
\ref{thm:semilocal} on certain complete semi-local algebras.
In the non-commutative case, it turns out that $M$ and $N$ have to be
$A$-bimodules, and symmetry of $\Ext$ vanishing takes the form
\[
  \Ext_A^i(M,N) = 0 \; \mbox{for} \; i \gg 0
  \; \Leftrightarrow \;
  \Ext_A^i(N,{}^{\sigma}\!M) = 0 \; \mbox{for} \; i \gg 0
\]
where $\sigma$ is a ``symmetrizing automorphism'' of $A$.

Note that while the theorems of this paper are phrased using only
classical homological algebra, some other parts of the paper use
derived categories and functors such as $\RHom$ and $\LTensor$.
However, the notation remains standard, and only standard properties
of derived categories are used.  Some background can be found in
\cite[sec.\ 2]{WZPI}, and that paper also explains the ring theory
notation which will be used.

The paper is organized as follows: After this introduction comes
section \ref{sec:main} which introduces properties (R) and (S) and
proves the abstract $\Ext$ vanishing result theorem \ref{thm:main}.
Next, sections \ref{sec:R} and \ref{sec:S} give some methods by which
properties (R) and (S) can be established.  And finally, section
\ref{sec:examples} uses the machinery to prove the concrete symmetry
theorems for $\Ext$ vanishing \ref{thm:commutative},
\ref{thm:Frobenius}, and \ref{thm:semilocal}.

\medskip
\noindent
{\bf Acknowledgement.}  
This paper was triggered by \cite{HJ} which inspired my attack on
symmetry of $\Ext$ vanishing by purely homological me\-thods.

The present version was rewritten completely after Izuru Mori alerted me
to an error in the original version.  I would like to thank professor
Mori for his help, and for communicating his preprint \cite{Mori}.

I thank Henrik Holm for a conversation related to lemma
\ref{lem:isomorphism}, and James Zhang for answering my questions
about property (S) over complete semi-local algebras.

\medskip

\section{An abstract $\Ext$ vanishing result}
\label{sec:main}

\begin{Definition}[Property (R)]
\label{dfn:R}
Let $A$ be a ring and let $N$ be an $A$-left-module.

Suppose that if $L$ is an exact complex of finitely generated
projective $A$-right-modules for which the cohomology of $L \otimes_A
N$ vanishes in high degrees, that is,
\[
  \H^{\gg 0}(L \otimes_A N) = 0,
\]
then in fact, all the cohomology of $L \otimes_A N$ vanishes, that is,
\[
  \H(L \otimes_A N) = 0.
\]

Then $N$ is said to have property (R).
\end{Definition}

\begin{Remark}
It is easy to see that if $N$ has finite flat dimension, then it has
property (R).  But there are other modules with property (R), see
section \ref{sec:R}.
\end{Remark}

\begin{Lemma}
\label{lem:isomorphism}
Let $A$ be a noetherian ring with finite injective dimension from the
left, $\id_A(A) < \infty$.  Let $M$ be a finitely generated
$A$-left-module and let $N$ be an $A$-bimodule.

Suppose that $N$ viewed as an $A$-left-module has property {\rm (R)}
and that $\Ext_A^{\gg 0}(M,N) = 0$.  Then
\[
  \RHom_A(M,A \LTensor_A N) 
  \cong 
  \RHom_A(M,A) \LTensor_A N
\]
in $\D(A^{\opp})$, the derived category of $A$-right-modules.
\end{Lemma}

\begin{proof}
Let
\[
  P = \cdots \rightarrow P^{-2} \rightarrow P^{-1} \rightarrow P^0
      \rightarrow 0 \rightarrow \cdots
\]
be a projective resolution of $M$ consisting of finitely generated
mo\-du\-les.  Then $P$ is a right-bounded complex of finitely generated
projective $A$-left-mo\-du\-les.  Hence $\Hom_A(P,A)$ is a
left-bounded complex of finitely generated projective
$A$-right-modules, and its cohomology is $\Ext_A(M,A)$ which is
bounded because $\id_A(A) < \infty$.  In consequence, there is a
quasi-isomorphism $Q \stackrel{\simeq}{\longrightarrow}
\Hom_A(P,A)$ where $Q$ is a right-bounded complex of finitely
generated projective $A$-right-modules.

The quasi-isomorphism can be completed to a distinguished triangle
\[
  Q \stackrel{\simeq}{\longrightarrow} \Hom_A(P,A) 
    \longrightarrow L \longrightarrow,
\]
where $L$ is now an exact complex of finitely generated projective
$A$-right-modules.  This again gives a distinguished triangle
\begin{equation}
\label{equ:z}
  Q \otimes_A N \longrightarrow \Hom_A(P,A) \otimes_A N
                \longrightarrow L \otimes_A N \longrightarrow.
\end{equation}

Since $Q$ and $P$ are projective resolutions of $\Hom_A(P,A)$ and $M$,
it follows that
\[
  Q \otimes_A N 
  \cong \Hom_A(P,A) \LTensor_A N
  \cong \RHom_A(M,A) \LTensor_A N
\]
and
\[
  \Hom_A(P,A) \otimes_A N
  \stackrel{\rm (a)}{\cong} \Hom_A(P,A \otimes_A N)
  \cong \RHom_A(M,A \LTensor_A N),
\]
where (a) is because each module in $P$ is finitely generated
projective.  So the distinguished triangle \eqref{equ:z} reads
\begin{equation}
\label{equ:a}
  \RHom_A(M,A) \LTensor_A N
  \longrightarrow \RHom_A(M,A \LTensor_A N)
  \longrightarrow L \otimes_A N \longrightarrow.
\end{equation}

The cohomology of $\RHom_A(M,A)$ is $\Ext_A(M,A)$ which is bounded
because $\id_A(A) < \infty$, as remarked above.  Consequently, the
cohomology of $\RHom_A(M,A) \LTensor_A N$ is right-bounded.  And the
cohomology of
\[
  \RHom_A(M,A \LTensor_A N) \cong \RHom_A(M,N)
\]
is $\Ext_A(M,N)$ which is right-bounded by assumption.

The distinguished triangle \eqref{equ:a} hence shows that the
cohomology of $L \otimes_A N$ is right-bounded, that is, $\H^{\gg 0}(L
\otimes_A N) = 0$.  Since $N$ viewed as an $A$-left-module has
property (R), it follows that $\H(L \otimes_A N) = 0$, and hence the
first morphism in \eqref{equ:a},
\[
  \RHom_A(M,A) \LTensor_A N
  \longrightarrow \RHom_A(M,A \LTensor_A N),
\]
is a quasi-isomorphism and hence an isomorphism in the derived
category $\D(A^{\opp})$.
\end{proof}

\begin{Remark}
Recall for the following setup that if $\sigma$ is an
automorphism of a ring $A$ and $M$ is an $A$-left-module, then there
is an $A$-left-module ${}^{\sigma}\!M$ with scalar multiplication
defined in terms of the scalar multiplication of $M$ by $a \cdot m =
\sigma(a)m$.  

This procedure can also be applied to $A$-right-modules, and from
either side to $A$-bimodules.  Note that as $A$-bimodules,
\[
  {}^{\sigma}\!A \cong A^{\sigma^{-1}}.
\]
\end{Remark}

\begin{Setup}
\label{set:Asigma}
Let $A$ be a ring with an automorphism $\sigma$ and suppose that the
$A$-bimodule ${}^{\sigma}\!A$ has a resolution
\[
  I = \cdots \rightarrow 0 \rightarrow I^0 \rightarrow I^1 \rightarrow
      I^2 \rightarrow \cdots
\]
of $A$-bimodules which are injective when viewed either as
$A$-left-mo\-du\-les or $A$-right-modules.
\end{Setup}

\begin{Remark}
If $I$ is viewed as a complex of $A$-left-modules, then it is an
injective resolution of $A$ viewed as an $A$-left-module.  Similarly
from the right.  Hence the functors
\[
  \RHom_A(-,{}^{\sigma}\!A) 
  \;\; \mbox{and} \;\;
  \RHom_{A^{\opp}}(-,{}^{\sigma}\!A)
\]
can be defined as 
\[
  \Hom_A(-,I)
  \;\; \mbox{and} \;\; 
  \Hom_{A^{\opp}}(-,I).  
\]
This has the advantage of giving $\RHom$ functors which are also
defined on the derived category of $A$-bimodules, hence enabling the
next definition.
\end{Remark}

\begin{Definition}[Property (S)]
\label{dfn:S}
Let $A$ and $\sigma$ be as in setup \ref{set:Asigma} and let $N$ be
an $A$-bimodule which is finitely generated from either side.  

Suppose that
\[
  \RHom_{A}(N,{}^{\sigma}\!A) \cong \RHom_{A^{\opp}}(N,{}^{\sigma}\!A)
\]
in $\D(A^{\opp})$.  Then $N$ is said to have property (S).
\end{Definition}

\begin{Remark}
It is clear that if $A$ is commutative and $\sigma$ is the identity,
then an $A$-module $N$, viewed as an $A$-bimodule via $an = na$ for
$a$ in $A$ and $n$ in $N$, has property (S).  But there are other
modules with property (S), see section \ref{sec:S}.
\end{Remark}

\begin{Theorem}
\label{thm:main}
Let $A$ and $\sigma$ be as in setup \ref{set:Asigma} and suppose that
$A$ is noetherian and that $\id_A(A) < \infty$ and $\id_{A^{\opp}}(A)
< \infty$.  Let $M$ be a finitely generated $A$-left-module and let
$N$ be an $A$-bimodule which is finitely generated from either side.

Suppose that $N$ viewed as an $A$-left-module has property {\rm (R)}
and that $N$ has property {\rm (S)}.  Then
\[
  \Ext_A^{\gg 0}(M,N) = 0
  \; \Rightarrow \;
  \Ext_A^{\gg 0}(N,{}^{\sigma}\!M) = 0.
\]
\end{Theorem}

\begin{proof}
I can compute,
\begin{align*}
  \lefteqn{
    \RHom_{A^{\opp}}(\RHom_A(M,N),{}^{\sigma}\!A)
          }
  & \\
  & \;\;\;\;\; \cong
    \RHom_{A^{\opp}}(\RHom_A(M,A \LTensor_A N),{}^{\sigma}\!A) \\
  & \;\;\;\;\; \stackrel{\rm (r)}{\cong}
    \RHom_{A^{\opp}}(\RHom_A(M,A) \LTensor_A N,{}^{\sigma}\!A) \\
  & \;\;\;\;\; \cong
    \RHom_{A^{\opp}}(\RHom_A(M,A),
                     \RHom_{A^{\opp}}(N,{}^{\sigma}\!A)) \\
  & \;\;\;\;\; \stackrel{\rm (s)}{\cong}
    \RHom_{A^{\opp}}(\RHom_A(M,A),
                     \RHom_A(N,{}^{\sigma}\!A)) \\
  & \;\;\;\;\; \cong
    \RHom_A(N,\RHom_{A^{\opp}}(\RHom_A(M,A),{}^{\sigma}\!A)) \\
  & \;\;\;\;\; \stackrel{\rm (t)}{\cong}
    \RHom_A(N,\RHom_{A^{\opp}}(A,{}^{\sigma}\!A) \LTensor_A M) \\
  & \;\;\;\;\; \cong
    \RHom_A(N,{}^{\sigma}\!A \LTensor_A M) \\
  & \;\;\;\;\; \cong
    \RHom_A(N,{}^{\sigma}\!M).
\end{align*}
Here (r) is by lemma \ref{lem:isomorphism} since $N$ viewed as an
$A$-left-module has pro\-per\-ty (R), and (s) is since $N$ has
property (S), while (t) is because ${}^{\sigma}\!A$ viewed as an
$A$-right-module is just $A$, so $\id_{A^{\opp}}({}^{\sigma}\!A) =
\id_{A^{\opp}}(A) < \infty$, and this implies
\[
  \RHom_{A^{\opp}}(\RHom_A(M,A),{}^{\sigma}\!A)
  \cong \RHom_{A^{\opp}}(A,{}^{\sigma}\!A) \LTensor_A M.
\]
The remaining isomorphisms are standard.

Now, the condition
\[
  \Ext_A^{\gg 0}(M,N) = 0
\]
says that the cohomology of
\[
  \RHom_A(M,N)
\]
is bounded.  But since $\id_{A^{\opp}}({}^{\sigma}\!A) < \infty$, this
implies that the cohomology of
\[
  \RHom_{A^{\opp}}(\RHom_A(M,N),{}^{\sigma}\!A)
\]
is bounded.  And then the above computation shows that the cohomology
of 
\[
  \RHom_A(N,{}^{\sigma}\!M)
\]
is bounded, that is,
\[
  \Ext_A^{\gg 0}(N,{}^{\sigma}\!M) = 0.
\]
\end{proof}

\begin{Remark}
Definition \ref{dfn:S} and theorem \ref{thm:main} would be simpler
without the automorphism $\sigma$.  However, it appears that to get a
theory covering any reasonable stock of non-commutative rings,
$\sigma$ is a necessary ingredient.  See remark \ref{rmk:Frobenius}
and proposition \ref{pro:semilocal}.
\end{Remark}

\section{Property {\rm (R)}}
\label{sec:R}

This section gives some methods by which property (R) from definition
\ref{dfn:R} can be established.

The following definition is classical.

\begin{Definition}
\label{dfn:rigid}
Let $A$ be a ring and let $N$ be an $A$-left-module.

Suppose that there exists $c \geqslant 0$ so that
\[
  \Tor^A_i(Z,N) = \cdots = \Tor^A_{i+c}(Z,N) = 0
  \; \Rightarrow \;
  \Tor^A_{\geqslant i}(Z,N) = 0
\]
for each $A$-right-module $Z$ and each $i \geqslant 1$.  Then $N$ is
said to be $\Tor$ rigid.
\end{Definition}

\begin{Proposition}
\label{pro:rigid}
Let $A$ be a ring and let $N$ be an $A$-left-module which is $\Tor$
rigid.

Then $N$ has property {\rm (R)}.
\end{Proposition}

\begin{proof}
Let $L$ be an exact complex of finitely generated projective
$A$-right-modules with 
\[
  \H^{\gg 0}(L \otimes_A N) = 0.
\]
Suspending $L$ if necessary, I can suppose
\begin{equation}
\label{equ:b}
  \H^{\geqslant -c-1}(L \otimes_A N) = 0
\end{equation}
where $c$ is the constant from definition \ref{dfn:rigid}.  

Since $L$ is exact,
\[
  \cdots \rightarrow L^{-2} 
  \rightarrow L^{-1}
  \rightarrow L^0
  \rightarrow Z 
  \rightarrow 0
\]
is a projective resolution of the $A$-right-module $Z =
\Coker\,(L^{-1} \rightarrow L^0)$, so
\begin{equation}
\label{equ:c}
  \Tor^A_j(Z,N) \cong \H^{-j}(L \otimes_A N)
\end{equation}
for $j \geqslant 1$.

Combining equations \eqref{equ:b} and \eqref{equ:c} shows
\[
  \Tor^A_1(Z,N) = \cdots = \Tor^A_{c+1}(Z,N) = 0,
\]
and $\Tor$ rigidity of $N$ now implies
\[
  \Tor^A_j(Z,N) = 0
\]
for each $j \geqslant 1$.  Combining with equation \eqref{equ:c} shows
\[
  \H^{\leqslant -1}(L \otimes_A N) = 0,
\]
and combining with equation \eqref{equ:b} shows
\[
  \H(L \otimes_A N) = 0
\]
as desired.
\end{proof}

The following definition was first made in the commutative case in
\cite[dfn.\ 3.1]{HJ}.

\begin{Definition}
[the AB property]
\label{dfn:AB}
Let $A$ be a ring.  Suppose that there exists $c \geqslant 1$ so that
\[
  \Ext_A^{\gg 0}(M,N) = 0
  \; \Rightarrow \;
  \Ext_A^{\geqslant c}(M,N) = 0
\]
when $M$ and $N$ are finitely generated $A$-left-modules.

Then $A$ is said to have the left AB property.
\end{Definition}

\begin{Proposition}
\label{pro:AB}
Let $A$ be a ring with $\id_{A^{\opp}}(A) < \infty$ which has the left
{\rm AB} property.

Then each finitely generated $A$-left-module $N$ has property {\rm
(R)}. 
\end{Proposition}

\begin{proof}
The condition $\id_{A^{\opp}}(A) < \infty$ implies that $A$ viewed as
an $A$-right-module has a bounded injective resolution,
\[
  I = \cdots \rightarrow 0 \rightarrow I^0 \rightarrow \cdots
  \rightarrow I^d \rightarrow 0 \rightarrow \cdots.
\]
There is a qua\-si-i\-so\-mor\-phism $A_A
\stackrel{\simeq}{\rightarrow} I$ where $A_A$ is $A$ viewed as an
$A$-right-module.  Let $L$ be an exact complex of finitely generated
projective $A$-right-modules.  Since $L$ consists of projective
modules, the functor $\Hom_{A^{\opp}}(L,-)$ preserves
quasi-isomorphisms of bounded complexes, so there is a
quasi-isomorphism
\[
  \Hom_{A^{\opp}}(L,A)
  \stackrel{\simeq}{\longrightarrow}
  \Hom_{A^{\opp}}(L,I).
\]
The complex $I$ is bounded and consists of injective modules, so
the functor $\Hom_{A^{\opp}}(-,I)$ preserves exactness, so
$\Hom_{A^{\opp}}(L,I)$ is exact.  Hence 
\[
  L^* = \Hom_{A^{\opp}}(L,A)
\]
is also exact.

Since $L$ consists of finitely generated projective $A$-right-modules,
there is an isomorphism $L \stackrel{\cong}{\longrightarrow}
\Hom_A(L^*,A)$.  Hence, if $N$ is a finitely generated
$A$-left-module,
\begin{align}
\nonumber
  L \otimes_A N 
  & \cong \Hom_A(L^*,A) \otimes_A N \\
\nonumber
  & \stackrel{\rm (a)}{\cong} \Hom_A(L^*,A \otimes_A N) \\
\label{equ:y}
  & \cong \Hom_A(L^*,N),
\end{align}
where (a) is because $L^*$ consists of finitely generated projective
$A$-left-modules.

Since $L^*$ is exact, for each $i$,
\[
  \cdots \rightarrow L^{* i - 2}
  \rightarrow L^{* i - 1}
  \rightarrow L^{* i}
  \rightarrow Z
  \rightarrow 0
\]
is a projective resolution of the finitely generated $A$-left-module
$Z = \Coker\,(L^{* i - 1} \rightarrow L^{* i})$.  Combining with
equation \eqref{equ:y} shows 
\begin{equation}
\label{equ:h}
  \H^j(L \otimes_A N)
  \cong \H^j \Hom_A(L^*,N)
  \cong \Ext_A^{j + i}(Z,N)
\end{equation}
for $j + i \geqslant 1$.  

Now suppose
\[
  \H^{\gg 0}(L \otimes_A N) = 0.  
\]
Then equation \eqref{equ:h} shows $\Ext_A^{\gg 0}(Z,N) = 0$, and the
left AB property implies $\Ext_A^{\geqslant c}(Z,N) = 0$ which by
equation \eqref{equ:h} says
\[
  \H^j(L \otimes_A N) = 0
\]
for $j + i \geqslant c$, that is, $j \geqslant c - i$.  Varying $i$
now shows
\[
  \H(L \otimes_A N) = 0
\]
as desired.
\end{proof}

\begin{Proposition}
\label{pro:CMtype}
Let $A$ be a left-noetherian ring with $\id_A(A) < \infty$ which has a
finite set $\cC$ of left-modules so that if $N$ is a finitely
generated $A$-left-module with $\Ext_A^{\geqslant 1}(N,A) = 0$, then $N$ is
isomorphic to the direct sum of finitely many modules from $\cC$.

Then $A$ has the left {\rm AB} property.
\end{Proposition}

\begin{proof}
The set $\cC$ is finite, so it is clear that there exists
$\widetilde{c} \geqslant 1$ so that if $\widetilde{M}$ and
$\widetilde{N}$ are isomorphic to direct sums of finitely many modules
from $\cC$, then
\begin{equation}
\label{equ:l}
  \Ext_A^{\gg 0}(\widetilde{M},\widetilde{N}) = 0
  \; \Rightarrow \;
  \Ext_A^{\geqslant \widetilde{c}}(\widetilde{M},\widetilde{N}) = 0.
\end{equation}

Let $M$ be a finitely generated $A$-left-module and write $d =
\id_A(A)$.  Let
\[
  0 \rightarrow \widetilde{M}
  \rightarrow P_{d-1}
  \rightarrow \cdots
  \rightarrow P_0
  \rightarrow M
  \rightarrow 0
\]
be an exact sequence where the $P_i$ are finitely generated projective
$A$-left-modules.  Since $\widetilde{M}$ is the $d$'th syzygy in a
projective resolution of $M$,
\begin{equation}
\label{equ:j}
  \Ext_A^{d+i}(M,N) \cong \Ext_A^i(\widetilde{M},N)
\end{equation}
for each $A$-left-module $N$ and each $i \geqslant 1$.  If $N = Q$ is
finitely generated projective, then this implies
\begin{equation}
\label{equ:k}
  \Ext_A^i(\widetilde{M},Q) = 0
\end{equation}
for $i \geqslant 1$, since $\id_A(Q) \leqslant \id_A(A) = d$.  In
particular, $\Ext_A^{\geqslant 1}(\widetilde{M},A) = 0$, so
$\widetilde{M}$ is isomorphic to the direct sum of finitely many
modules from $\cC$. 

Now let $N$ be any finitely generated $A$-left-module, and let
\[
  0 \rightarrow \widetilde{N}
  \rightarrow Q_{d-1}
  \rightarrow \cdots
  \rightarrow Q_0
  \rightarrow N
  \rightarrow 0
\]
be an exact sequence where the $Q_i$ are finitely generated projective
$A$-left-modules.  Like $\widetilde{M}$, the module $\widetilde{N}$ is
isomorphic to the direct sum of finitely many modules from $\cC$.
Splitting the exact sequence into short exact sequences and applying
the long exact sequence of $\Ext$ groups repeatedly along with
equation \eqref{equ:k} shows
\[
  \Ext_A^i(\widetilde{M},N) 
  \cong \Ext_A^{d+i}(\widetilde{M},\widetilde{N})
\]
for $i \geqslant 1$, and combining this with equation \eqref{equ:j} shows
\begin{equation}
\label{equ:m}
  \Ext_A^j(M,N) \cong \Ext_A^j(\widetilde{M},\widetilde{N})
  \;\; \mbox{for} \;\; j \geqslant d+1.
\end{equation}

Setting $c = \max \{\, \widetilde{c},d+1 \,\}$ and combining equations
\eqref{equ:l} and \eqref{equ:m} now shows 
\[
  \Ext_A^{\gg 0}(M,N) = 0
  \; \Rightarrow \;
  \Ext_A^{\geqslant c}(M,N) = 0
\]
as desired.
\end{proof}

\section{Property {\rm (S)}}
\label{sec:S}

This section gives some methods by which property (S) from definition
\ref{dfn:S} can be established.

\begin{Remark}
\label{rmk:commutative}
Let $A$ be a commutative ring and let $\sigma$ be the identity
automorphism of $A$.  This clearly gives the situation of setup
\ref{set:Asigma}.

If $N$ is a finitely generated $A$-module, then I can view $N$ as an
$A$-bimodule via $an = na$ for $a$ in $A$ and $n$ in $N$,
and it is obvious that
\[
  \RHom_A(N,A) \cong \RHom_{A^{\opp}}(N,A)
\]
in $\D(A^{\opp})$, so $N$ has property (S).
\end{Remark}

\begin{Remark}
\label{rmk:Frobenius}
Let $A$ be a (finite dimensional) Frobenius algebra over the field $k$
and let $\sigma$ be a symmetrizing automorphism of $A$.  This means
that $\sigma$ is an automorphism for which $\Hom_k(A,k) \cong
{}^{\sigma}\!A$ as $A$-bimodules.

Observe that ${}^{\sigma}\!A$ is injective as an $A$-module from
either side, so ${}^{\sigma}\!A$ is a resolution of itself by
$A$-bimodules which are injective when viewed either as
$A$-left-modules or as $A$-right-modules.  Hence I am in the situation
of setup \ref{set:Asigma}.

Let $N$ be an $A$-bimodule which is finitely generated from either
side.  Then
\begin{align*}
  \RHom_A(N,{}^{\sigma}\!A) 
  & \cong \RHom_A(N,\Hom_k(A,k)) \\
  & \cong \Hom_k(A \LTensor_A N,k) \\
  & \cong \Hom_k(N \LTensor_A A,k) \\
  & \cong \RHom_{A^{\opp}}(N,\Hom_k(A,k)) \\
  & \cong \RHom_{A^{\opp}}(N,{}^{\sigma}\!A)
\end{align*}
in $\D(A^{\opp})$, so $N$ has property (S).
\end{Remark}

\begin{Setup}
\label{set:semilocal}
Let $k$ be a field and let $A$ be a complete semi-local noetherian
$k$-algebra.
\end{Setup}

\begin{Remark}
\label{rmk:Morita}
The conditions of completeness and semi-locality in the setup mean
that, if $\fm$ denotes the Jacobson radical of $A$, then $A$ is
complete in the $\fm$-adic topology while $A/\fm$ is semi-simple with
$\dim_k A/\fm < \infty$.

Duality with respect to $k$ will be denoted by
\[
  (-)^{\prime} = \Hom_k(-,k).  
\]
This functor interchanges $A$-left-modules and $A$-right-modules, and
sends $A$-bimodules to $A$-bimodules.

According to \cite[lem.\ 2.5]{CWZ}, the $A$-bimodule
\[
  E = \colim (A/\fm^n)^{\prime}
\]
is injective from either side and induces a Morita self-duality for
$A$; see \cite[\S 24]{AndersonFuller} for background on Morita
duality.  

The Morita duality functors will be denoted by
\[
  (-)^L = \Hom_A(-,E),
  \;\;\;
  (-)^R = \Hom_{A^{\opp}}(-,E).
\]
The functor $(-)^L$ sends $A$-left-modules to $A$-right-modules and
$A$-bimodules to $A$-bimodules.  It sends modules finitely generated
from the left to modules artinian from the right, and modules artinian
from the left to modules finitely generated from the right; see
\cite[prop.\ 10.10 and thms.\ 24.5 and 24.6]{AndersonFuller}.  If $X$
is an $A$-left-module or an $A$-bimodule and $X$ is either finitely
generated or artinian from the left, then there is an isomorphism
\begin{equation}
\label{equ:i}
  (X^L)^R \cong X
\end{equation}
by \cite[prop.\ 10.10 and thm.\ 24.6]{AndersonFuller}.  Of course,
this can all be dualized.

The functors $(-)^{\prime}$, $(-)^L$, and $(-)^R$ are exact, and so
remain well defined on derived categories.
\end{Remark}

\begin{Lemma}
\label{lem:Lprime}
Let $A$ be as in setup \ref{set:semilocal} and let $Y$ be an
$A$-bimodule which is finitely generated from the right.

Then there is an isomorphism of $A$-bimodules
\[
  (Y^R)^{\prime} \cong Y.
\]
\end{Lemma}

\begin{proof}
First observe
\begin{align*}
  Y^R
  & = \Hom_{A^{\opp}}(Y,E) \\
  & = \Hom_{A^{\opp}}(Y,\colim (A/\fm^n)^{\prime}) \\
  & \stackrel{\rm (a)}{\cong} \colim \Hom_{A^{\opp}}(Y,(A/\fm^n)^{\prime}) \\
  & = \colim \Hom_{A^{\opp}}(Y,\Hom_k(A/\fm^n,k)) \\
  & \cong \colim \Hom_k(Y \otimes_A A/\fm^n,k) \\
  & = \colim\,(Y \otimes_A A/\fm^n)^{\prime} \\
  & \cong \colim\,(Y/Y\fm^n)^{\prime},
\end{align*}
where (a) is because $Y$ is finitely generated from the right
(this was remarked already in \cite[lem.\ 5.9(1)]{CWZ}).
Observe also that since $\dim_k A/\fm < \infty$, it is more generally
true that $\dim_k A/\fm^n < \infty$ for each $n \geqslant 1$.  Since
$Y/Y\fm^n$ is finitely generated from the right over $A/\fm^n$, it
follows that $\dim_k Y/Y\fm^n < \infty$ for each $n \geqslant 1$.

These observations imply (b) and (c) in
\begin{align*}
  (Y^R)^{\prime}
  & = \Hom_k(Y^R,k) \\
  & \stackrel{\rm (b)}{\cong} \Hom_k(\colim\,(Y/Y\fm^n)^{\prime},k) \\
  & \cong \lim \Hom_k((Y/Y\fm^n)^{\prime},k) \\
  & \stackrel{\rm (c)}{\cong} \lim Y/Y\fm^n \\
  & = \widehat{Y} \\
  & \stackrel{\rm (d)}{\cong} Y \otimes_A \widehat{A} \\
  & \stackrel{\rm (e)}{\cong} Y,
\end{align*}
where $\widehat{Y}$ denotes the completion of $Y$ in the $\fm$-adic
topology from the right, where (d) is by \cite[lem.\ 2.4(1)]{WZPI},
and where (e) holds because $\widehat{A} \cong A$ by assumption.
\end{proof}

\begin{Lemma}
\label{lem:LR}
Let $A$ be as in setup \ref{set:semilocal} and let $Z$ be an
$A$-bimodule which is artinian from either side.

Then there is an isomorphism of $A$-bimodules
\[
  Z^L \cong Z^R.
\]
\end{Lemma}

\begin{proof}
This is a computation,
\[
  Z^L
  \stackrel{\rm (a)}{\cong} ((Z^L)^R)^{\prime}
  \stackrel{\rm (b)}{\cong} Z^{\prime}
  \stackrel{\rm (c)}{\cong} ((Z^R)^L)^{\prime}
  \stackrel{\rm (d)}{\cong} Z^R.
\]
Here (a) holds by lemma \ref{lem:Lprime} because $Z^L$ is finitely
generated from the right by remark \ref{rmk:Morita}, and (b) is by
equation \eqref{equ:i} because $Z$ is artinian from the left.
Similarly, (c) is by the dual of  \eqref{equ:i} because $Z$ is
artinian from the right, and (d) is by the dual of lemma
\ref{lem:Lprime} because $Z^R$ is finitely generated from the left by
remark \ref{rmk:Morita}.
\end{proof}

\begin{Remark}
\label{rmk:Gammam}
The $\fm$-left-torsion functor over $A$ is defined by
\[
  \Gammam(M) = \{\, m \in M \,|\, \fm^n m = 0 
               \; \mbox{for} \; n \gg 0 \,\}.
\]
This functor sends $A$-left-modules to $A$-left-modules and
$A$-bimodules to $A$-bimodules.  It has a derived functor $\RGammam$
defined on derived ca\-te\-go\-ri\-es.  

There is also an $\fm$-right-torsion functor $\Gammamopp$ with derived
functor $\RGammamopp$.  See \cite[sec.\ 6]{WZPI}.
\end{Remark}

\begin{Lemma}
\label{lem:sym}
Let $A$ be as in setup \ref{set:semilocal}, suppose that $A$ has a
balanced dualizing complex (see \cite[dfn.\ 3.7]{CWZ}), and let $N$ be
an $A$-bimodule which is finitely generated from either side.

Then 
\[
  \RHom_A(N,\RGammam(A)^L) \cong \RHom_{A^{\opp}}(N,\RGammamopp(A)^R)
\]
in $\D(A \otimes_k A^{\opp})$, the derived category of $A$-bimodules. 
\end{Lemma}

\begin{proof}
The cohomology modules of $\RGammam(N)$ are artinian from either side
by the proof of \cite[thm.\ 3.5(3)]{CWZ}, and so lemma \ref{lem:LR}
implies
\begin{equation}
\label{equ:x}
  \RGammam(N)^L \cong \RGammam(N)^R
\end{equation}
in $\D(A \otimes_k A^{\opp})$. 

Since $A$ has a balanced dualizing complex, it satisfies the $\chi$
condition, as defined in \cite[p.\ 289]{CWZ}, by \cite[cor.\
3.9]{CWZ}.  Also, $A$ is complete in the $\fm$-adic topology, so $\fm$
has the left and right Artin-Rees properties by \cite[thm.\ 1.1]{Ja}.
Hence
\[
  \RGammam(N) \cong \RGammamopp(N)
\]
in $\D(A \otimes_k A^{\opp})$ by \cite[thm.\ 2.9]{WZDual}, so
\[
  \RGammam(N)^R \cong \RGammamopp(N)^R.
\]

Combining with equation \eqref{equ:x} gives
\[
  \RGammam(N)^L \cong \RGammamopp(N)^R,
\]
yielding (b) in
\begin{align*}
  \RHom_A(N,\RGammam(A)^L)
  & \stackrel{\rm (a)}{\cong} \RGammam(N)^L \\
  & \stackrel{\rm (b)}{\cong} \RGammamopp(N)^R \\
  & \stackrel{\rm (c)}{\cong} \RHom_{A^{\opp}}(N,\RGammamopp(A)^R)
\end{align*}
where (a) and (c) are by the version of local duality in \cite[thm.\
3.6(2)]{WZDual}.
\end{proof}

\begin{Proposition}
\label{pro:semilocal}
Let $A$ be a complete semi-local noetherian $k$-algebra over the field
$k$, which is Gorenstein in the sense that there is an automorphism
$\sigma$ so that the $d$'th suspension $\Sigma^d({}^{\sigma}\!A)$ is a
balanced dualizing complex (see \cite[dfn.\ 3.7]{CWZ}).  

Then $A$ and $\sigma$ give the situation of setup
\ref{set:Asigma}, and each $A$-bimodule $N$ which is finitely
generated from either side has property {\rm (S)}.
\end{Proposition}

\begin{proof}
To get the resolution $I$ in setup \ref{set:Asigma}, view the
$A$-bimodule ${}^{\sigma}\!A$ as a left-module over the enveloping
algebra $A \otimes_k A^{\opp}$, take an injective resolution $I$, and
view $I$ as a complex of $A$-bimodules.  


By the proof of \cite[cor.\ 3.9]{CWZ}, both $\RGammam(A)^L$ and
$\RGammamopp(A)^R$ are ba\-lan\-ced dualizing complexes for $A$.  But
the balanced dualizing complex is unique by \cite[p.\ 300]{CWZ}, so
\[
  \RGammam(A)^L \cong \Sigma^d({}^{\sigma}\!A) \cong \RGammamopp(A)^R.
\]
Substituting this into lemma \ref{lem:sym} gives
\[
  \RHom_A(N,\Sigma^d({}^{\sigma}\!A))
  \cong \RHom_{A^{\opp}}(N,\Sigma^d({}^{\sigma}\!A))
\]
in $\D(A \otimes_k A^{\opp})$, hence in particular in $\D(A^{\opp})$.
Using $\Sigma^{-d}$ then proves property (S) for $N$.
\end{proof}

\section{Symmetry theorems for $\Ext$ vanishing}
\label{sec:examples}

This section applies the theory of the previous sections to show
symmetry theorems for $\Ext$ vanishing over commutative rings,
Frobenius algebras, and (non-commutative) complete semi-local
algebras.

The following theorem on commutative rings uses complete intersection
dimension, $\CIdim$, as introduced in \cite{AGP}.  Note that a module
over a commutative ring can be viewed either as a left-module, a
right-module, or a bimodule via $an = na$ for $a$ in $A$ and $n$ in
$N$.  In consequence, the left AB property will just be called the AB
property.

\begin{Theorem}
\label{thm:commutative}
Let $A$ be a commutative local noetherian Gorenstein ring, and let $M$
and $N$ be finitely generated $A$-modules.  

Suppose either that 
\[
  \CIdim_A M < \infty \; \mbox{and} \:
  \CIdim_A N < \infty,
\]
or that $A$ has the {\rm AB} property (see definition \ref{dfn:AB}).
Then 
\[
  \Ext_A^{\gg 0}(M,N) = 0
  \; \Leftrightarrow \;
  \Ext_A^{\gg 0}(N,M) = 0.
\]

Moreover, if $A$ has finite Cohen-Macaulay type, then it has the {\rm
AB} property. 
\end{Theorem}

\begin{proof}
The proof is an application of theorem \ref{thm:main}, so I need to
check the conditions of that theorem.  

Let $\sigma$ be the identity automorphism.  This clearly gives
the situation of setup \ref{set:Asigma}.

The ring $A$ is commutative, local, and Gorenstein, so $\id_A(A) <
\infty$ and $\id_{A^{\opp}}(A) < \infty$.

Remark \ref{rmk:commutative} says that $M$ and $N$, viewed as
$A$-bimodules, have property (S).

If $\CIdim_A M < \infty$ and $\CIdim_A N < \infty$, then \cite[cor.\
2.3]{Jo} implies that $M$ and $N$ are $\Tor$ rigid in the sense of
definition \ref{dfn:rigid}, and then $M$ and $N$, viewed as
$A$-left-modules, have property (R) by proposition \ref{pro:rigid}.
If $A$ has the AB property then $M$ and $N$, viewed as
$A$-left-modules, have property (R) by proposition \ref{pro:AB}. 

So theorem \ref{thm:main} gives
\[
  \Ext_A^{\gg 0}(M,N) = 0
  \; \Rightarrow \;
  \Ext_A^{\gg 0}(N,M) = 0,
\]
and theorem \ref{thm:main} applied to $N$ and $M$ gives
\[
  \Ext_A^{\gg 0}(N,M) = 0
  \; \Rightarrow \;
  \Ext_A^{\gg 0}(M,N) = 0.
\]
This establishes the implications of the theorem.

Finally, a finitely generated $A$-module $N$ with $\Ext_A^{\geqslant
1}(N,A) = 0$ is maximal Cohen-Macaulay by \cite[thm.\ 3.3.7 and cor.\
3.5.11]{BH}, so such an $N$ is a direct sum of finitely many
indecomposable maximal Cohen-Macaulay modules.  Thus, if $A$ has
finite Cohen-Macaulay type, the set $\cC$ in proposition
\ref{pro:CMtype} can be taken to be a set of representatives of the
finitely many isomorphism classes of indecomposable maximal
Cohen-Macaulay modules, and so in this case, $A$ has the AB property.
\end{proof}

\begin{Remark}
The AB case of theorem \ref{thm:commutative} is known already from
\cite{HJ}, but the $\CIdim$ and finite Cohen-Macaulay type cases
appear to be new.

Note that a local complete intersection ring has the AB property by
\cite[cor.\ 3.5]{HJ}, so there is a non-trivial supply of rings with
the AB property.  Note also that the $\CIdim$ of any finitely
generated module over a local complete intersection ring is finite by 
\cite[thm.\ 1.3]{AGP}, so again, there is a non-trivial supply of
modules with finite $\CIdim$.
\end{Remark}

\begin{Theorem}
\label{thm:Frobenius}
Let $A$ be a (finite dimensional) Frobenius algebra over the field
$k$, let $\sigma$ be a symmetrizing automorphism of $A$, that is, an
automorphism for which $\Hom_k(A,k) \cong {}^{\sigma}\!A$ as
$A$-bimodules, and let $M$ and $N$ be $A$-bimodules which are finitely
generated from either side.

Suppose that $A$ has the left {\rm AB} property.  Then
\[
  \Ext_A^{\gg 0}(M,N) = 0
  \; \Leftrightarrow \;
  \Ext_A^{\gg 0}(N,{}^{\sigma}\!M) = 0.
\]

Moreover, if $A$ has finite representation type, then it has the left
{\rm AB} property.
\end{Theorem}

\begin{proof}
The proof is again an application of theorem \ref{thm:main}.

The algebra $A$ and the automorphism $\sigma$ give the
situation of setup \ref{set:Asigma} by remark \ref{rmk:Frobenius}.

The algebra $A$ is Frobenius, so in particular self-injective from
either side, so $\id_A(A) = 0 < \infty$ and $\id_{A^{\opp}}(A) = 0 <
\infty$.

Since $A$ has the left AB property, $N$ and ${}^{\sigma}\!M$ viewed as
$A$-left-modules have property (R) by proposition \ref{pro:AB}. 

The $A$-bimodules $N$ and ${}^{\sigma}\!M$ have property (S)
by remark \ref{rmk:Frobenius}. 

So theorem \ref{thm:main} gives
\[
  \Ext_A^{\gg 0}(M,N) = 0
  \; \Rightarrow \;
  \Ext_A^{\gg 0}(N,{}^{\sigma}\!M) = 0,
\]
and theorem \ref{thm:main} applied to $N$ and ${}^{\sigma}\!M$
gives
\[
  \Ext_A^{\gg 0}(N,{}^{\sigma}\!M) = 0
  \; \Rightarrow \;
  \Ext_A^{\gg 0}({}^{\sigma}\!M,{}^{\sigma}\!N) = 0,
\]
that is,
\[
  \Ext_A^{\gg 0}(N,{}^{\sigma}\!M) = 0
  \; \Rightarrow \;
  \Ext_A^{\gg 0}(M,N) = 0.
\]
This establishes the implications of the theorem.

Finally, if $A$ has finite representation type, then proposition
\ref{pro:CMtype} clearly implies that $A$ has the AB property.
\end{proof}

\begin{Remark}
There is a significant supply of Frobenius algebras of finite
representation type to which theorem \ref{thm:Frobenius} applies, for
instance, group algebras $kG$ where $k$ is a field of characteristic
$p > 0$ and $G$ is a finite group whose order is divisible by $p$ and
whose Sylow $p$-subgroups are cyclic, cf.\ \cite[thm.\ VI.3.3]{ARS}.
Here the automorphism $\sigma$ is even the identity since $kG$ is a
symmetric algebra.
\end{Remark}

The following definition was basically made already in \cite[dfn.\
8.5]{Chan} and \cite[dfn.\ 9.9]{Chan}.

\begin{Definition}
Let $A$ be a complete semi-local noetherian $k$-algebra over the field
$k$ for which the $\fm$-left-torsion functor $\Gammam$ from remark
\ref{rmk:Gammam} has finite cohomological dimension $d$.

A finitely generated $A$-left-module $M$ is called maximal
Co\-hen-Ma\-cau\-lay if
\[
  \R^{\leqslant d-1} \Gammam(M) = 0.
\]
If there are only finitely many isomorphism classes of indecomposable
maximal Cohen-Macaulay $A$-left-modules, then $A$ is said to have
finite Cohen-Macaulay type.
\end{Definition}

\begin{Theorem}
\label{thm:semilocal}
Let $A$ be a complete semi-local noetherian algebra over the field
$k$, which is Gorenstein in the sense that there is an
automorphism $\sigma$ of $A$ so that the $d$'th suspension
$\Sigma^d({}^{\sigma}\!A)$ is a balanced dualizing complex (see
\cite[dfn.\ 3.7]{CWZ}), and let $M$ and $N$ be $A$-bimodules which are
finitely generated from either side.

Suppose that $A$ has the left {\rm AB} property.  Then
\[
  \Ext_A^{\gg 0}(M,N) = 0
  \; \Leftrightarrow \;
  \Ext_A^{\gg 0}(N,{}^{\sigma}\!M) = 0.
\]

Moreover, if $A$ has finite Cohen-Macaulay type, then it has the left
{\rm AB} property.
\end{Theorem}

\begin{proof}
The proof is once more an application of theorem \ref{thm:main}.

The algebra $A$ and the automorphism $\sigma$ give the
situation of setup \ref{set:Asigma} by proposition
\ref{pro:semilocal}. 

The balanced dualizing complex $\Sigma^d({}^{\sigma}\!A)$ has finite
injective dimension from either side, so the same holds for
${}^{\sigma}\!A$.  As an $A$-right-module, ${}^{\sigma}\!A$ is just $A$,
so this implies $\id_{A^{\opp}}(A) < \infty$.  As an $A$-bimodule,
${}^{\sigma}\!A$ is $A^{\sigma^{-1}}$, and as an $A$-left-module, this
is $A$, so it also implies $\id_A(A) < \infty$.

Since $A$ has the left AB property, $N$ and ${}^{\sigma}\!M$ viewed as
$A$-left-modules have property (R) by proposition
\ref{pro:AB}. 

The $A$-bimodules $N$ and ${}^{\sigma}\!M$ have property (S) by
proposition \ref{pro:semilocal}.

Hence theorem \ref{thm:main} establishes the implications of the
theorem, just as in the proof of theorem \ref{thm:Frobenius}.

Finally, the balanced dualizing complex $\Sigma^d({}^{\sigma}\!A)
\cong \Sigma^d(A^{\sigma^{-1}})$ is also pre-balanced by \cite[lem.\
3.3]{CWZ}, so \cite[prop.\ 3.4]{CWZ} gives
\[
  \RGammam(N)^L \cong \RHom_A(N,\Sigma^d(A^{\sigma^{-1}}))
\]
when $N$ is a finitely generated $A$-left-module.  The $(-i)$'th
cohomology of this is
\[
  \R^i\Gammam(N)^L \cong \Ext_A^{d-i}(N,A^{\sigma^{-1}})
                   \cong \Ext_A^{d-i}(N,A)^{\sigma^{-1}}.
\]
This formula makes it clear that
\[
  N \; \mbox{is maximal Cohen-Macaulay}
  \; \Leftrightarrow \;
  \Ext_A^{\geqslant 1}(N,A) = 0,
\]
and so, if $A$ has finite Cohen-Macaulay type, proposition
\ref{pro:CMtype} with $\cC$ equal to a set of representatives of the
isomorphism classes of indecomposable maximal Cohen-Macaulay modules
implies that $A$ has the left AB property.
\end{proof}

\begin{Remark}
There is a supply of algebras to which theorem \ref{thm:semilocal}
applies.  Namely, consider the non-commutative special quotient
surface singularities introduced in \cite{Chan}.  These are fixed
point algebras of the form $A = B^G$ where $G$ is a finite group
acting suitably on the (non-commutative) complete local noetherian
regular algebra $B$; see \cite[sec.\ 3]{Chan}.

Such an $A$ is complete local noetherian.  It follows from
\cite[cor.\ 6.11, dfn.\ 6.15, and prop.\ 9.8]{Chan} that
$\Sigma^2({}^{\sigma}\!A)$ is a balanced dualizing complex for $A$
for some automorphism $\sigma$.  And $A$ has finite Cohen-Macaulay
type by \cite[thm.\ 9.10]{Chan}, so theorem \ref{thm:semilocal}
applies. 

In another direction, \cite[cor.\ 2.4]{Mori} says that if $A$ is a
noetherian ring with finite global dimension, then $A/(x_1, \ldots,
x_n)$ has the left AB property when $x_1, \ldots, x_n$ is a regular
sequence of normal elements where conjugation by $x_{i+1}$ has finite
order on $A/(x_1, \ldots, x_i)$ for each $i$.  This gives another way
of getting algebras to which theorems \ref{thm:Frobenius} and
\ref{thm:semilocal} apply.
\end{Remark}


\end{document}